\documentclass[12pt,reqno]{amsart}

\usepackage{amssymb}
\usepackage{amscd}
\usepackage{amsfonts}
\usepackage{setspace}
\usepackage{version}

\usepackage{graphicx}

\newtheorem{theorem}{Theorem}

\theoremstyle{definition}

\newtheorem{conjecture}[theorem]{Conjecture}

\theoremstyle{definition}

\newcommand{\md}[1]{\ensuremath{(\operatorname{mod}\, #1)}}
\newcommand{\mdsub}[1]{\ensuremath{(\mbox{\scriptsize mod}\, #1)}}

\renewcommand{\leq}{\leqslant}
\renewcommand{\geq}{\geqslant}

\newcommand\weak{\operatorname{weak}}
\def\F{\mathbb{F}}

\def\Z{\mathbb{Z}}
\def\E{\mathbb{E}}

\newcommand\EH{\operatorname{EH}}
\def\eps{\varepsilon}

\begin{document}

\title[Bounded gaps between primes]{Bounded gaps between primes}


\author{Ben Green}
\address{Mathematical Institute, Andrew Wiles Building, Radcliffe Observatory Quarter, Woodstock Rd, Oxford OX2 6GG.}
\email{ben.green@maths.ox.ac.uk}

\onehalfspace

\begin{abstract}
These are notes on Zhang's work and subsequent developments produced in preparation for 5 hours of talks for a general mathematical audience given in Cambridge, Edinburgh and Auckland over the last year.
\end{abstract}

\maketitle
In May 2013, Yitang Zhang stunned the mathematical world by proving the following result, which we shall loosely refer to as ``bounded gaps between primes''.

\begin{theorem}[Zhang]
There is an absolute constant $H$ such that there are infinitely many pairs of distinct primes differing by at most $H$. 
\end{theorem}

The celebrated twin prime conjecture asserts that one may take $H = 2$. Zhang obtained $H = 70000000$, but subsequent feverish activity by a massively collaborative \emph{Polymath} project reduced this to 4680. In late 2013, James Maynard and Terry Tao found a much simpler proof of Zhang's result giving $H = 600$, and a further \emph{Polymath} project based on this work has, at the time of writing, reduced $H$ to 252.

These notes arose as a result of preparing for three general audience lectures on Zhang's work, in Cambridge and Edinburgh in June 2013 and Auckland in February 2014. On each occasion the lectures were for a general mathematical audience, so I skipped essentially every detail\footnote{In particular, in my attitude to $\log X$ I had in mind the remarks about formulae allegedly made to Stephen Hawking by the editor of \emph{A Brief History of Time}: each occurrence would halve his readership.} of the argument and said many things that were technically false, or at least rather vague, even if ``morally correct''. My main aim was to give as much of a flavour as possible of all the main ingredients in about 2 hours (or, in Edinburgh, just 1 hour). I also tried to mention all the key players in the story, but did not give precise citations or references.  I have retained this style here. The first \emph{Polymath} paper is, though long, very well-written and gives all of these details. It is freely available at 
\[ \mbox{\texttt{http://arxiv.org/abs/1402.0811}}.\] A highly recommended resource for a more detailed account of the result is Andrew Granville's article for the 2014 \emph{Current Events} session of the AMS, freely available at \[ \mbox{\texttt{http://www.dms.umontreal.ca/$\sim$andrew/CEBBrochureFinal.pdf}}.\] Many people, including me, benefitted greatly from Kannan Soundarajan's beautiful article on the work of Goldston, Pintz and Y{\i}ld{\i}r{\i}m for the 2006 edition of the same session, freely available at \[ \mbox{\texttt{http://arxiv.org/abs/math/0605696}.}\]

Something I wish to emphasise in these notes is that Zhang's result should be thought of as the culmination of ideas about prime numbers developed by many of the great analytic number theorists of the 20th century. It melds two important bodies of work:

1. Ideas of Goldston, Pintz and Y{\i}ld{\i}r{\i}m, building on work of Selberg and others, establishing a link between the distribution of primes in arithmetic progressions and small gaps between primes;

2. Ideas of Bombieri, Fouvry, Friedlander and Iwaniec, building on but going well beyond work of Bombieri-Vinogradov, pinning down strong results about how primes are distributed in arithmetic progressions.\vspace{11pt}

\textsc{Primes in progressions.} The distribution of primes in progressions is central to the whole story, so let us begin with a brief tour of that subject. We begin by recalling the prime number theorem, which states that $\pi(X)$, the number of primes less than or equal to $X$, is roughly $X/\log X$. The fact that $\log X \rightarrow \infty$, which means the primes have density tending to zero, of course explains why Zhang's theorem is not at all obvious: the \emph{average} gap between primes less than $X$ is about $\log X$. One might also remark that Zhang's theorem is not at all \emph{surprising}, either, since a random set of $X/\log X$ integers less than $X$ will, for $X$ large, have many pairs spaced by at most $2$. In fact, it will have many pairs spaced by at most $1$, something not true for the primes themselves of course -- to model the primes by a random set one has to be more careful and take account, for example, of the fact that most primes are odd.

It is convenient to state the prime number theorem a little differently by introducing the \emph{von Mangoldt function} $\Lambda(n)$, which is basically defined to equal $\log n$ if $n$ is a prime\footnote{And $\log p$ if $n = p^k$ is a prime power, but this detail will not bother us here.}. Then the prime number theorem is equivalent to saying that 
\[ \sum_{x \leq X} \Lambda(x) \approx X.\]
How many primes $x$ satisfy some additional congruence condition $x \equiv c \md{d}$? Clearly (apart from for very small primes $x \leq d$) we must have $c$ coprime to $d$, but there is no obvious additional restriction. In fact, one expects the primes to be equally distributed amongst the $\phi(d)$ residue classes coprime to $d$. If this is the case, we have
\[ \sum_{\substack{x \leq X \\ x \equiv c \mdsub{d}}} \Lambda(x) \approx \frac{X}{\phi(d)}.\]
When this is true for a given value of $d$ and for \emph{all} $c$ coprime to $d$ then we shall say that the primes are \emph{nicely distributed}\footnote{Of course, this is only an informal definition and not a rigorous one because we have not elaborated upon the meaning of $\approx$.} modulo $d$. 

Proving this statement is another matter. Think of $X$ tending to infinity: then one would like to know that the primes are nicely distributed modulo $d$ for all $d$ up to some limit, which it would be nice to take as large as possible. This is, however, only known when $d$ is less than a power of $\log X$ (in fact $d$ can be less than \emph{any} power of $\log X$, a statement known as the Siegel-Walfisz theorem). 

It turns out that proving that the primes are nicely distributed for $d$ up to about $X^{1/2}$ is equivalent to the Generalised Riemann Hypothesis, so one should not expect too much progress soon.

Remarkably, one can do far better if one is prepared to know that the primes are nicely distributed only for\footnote{Once again, this is an informal definition. To make it rigorous, one would need to elaborate upon the meaning of \emph{almost all} as well as the $\approx$ notation from earlier.} \emph{almost all} $d$. The classic result in this vein is the Bombieri-Vinogradov theorem, which asserts that the primes are nicely distributed modulo $d$ for almost all $d \leq X^{1/2}$. The Bombieri-Vinogradov theorem is often described as a kind of \emph{Riemann Hypothesis on Average}. 

It is suspected that even more is true, a conjecture known as the Elliott-Halberstam Conjecture. There is a different Elliott-Halberstam conjecture $\EH(\theta)$ for each value of $\theta < 1$, and these conjectures get stronger as $\theta$ increases. Here is a rough statement.
\begin{conjecture}[Elliott Halberstam conjecture $\EH(\theta)$]
The primes are nicely distributed modulo $d$ for most values of $d$ up to $X^{\theta}$.
\end{conjecture}
Currently, we cannot prove this statement for any value of $\theta > \frac{1}{2}$. \vspace{11pt}

\textsc{GPY and BFFI.} Goldston, Pintz and Y{\i}ld{\i}r{\i}m (henceforth referred to as GPY) established a remarkable link between the problem of finding bounded gaps between primes and the Elliott-Halberstam conjecture. 
\begin{theorem}[GPY, 2005]
Suppose the Elliott-Halberstam conjecture $\EH(\theta)$ holds for some value of $\theta > \frac{1}{2}$. That is, suppose the primes are nicely distributed modulo $d$ for most values of $d$ up to $X^{\theta}$. Then we have bounded gaps between primes. 
\end{theorem}
One should also say that GPY \emph{unconditionally} proved some results about gaps between primes far superior to any that had appeared before their work, showing for example that there are always pairs of primes of size around $X$ and separated by about $\sqrt{\log X}$, far less than the average spacing of $\log X$ for primes of this size. 

About 20 years before that, in the 1980s, deep work of Bombieri, Fouvry, Friedlander and Iwaniec (in various combinations) had shown that a certain weak variant of the Elliott-Halberstam conjecture is true for some $\theta > \frac{1}{2}$. In fact they obtained $\theta = \frac{4}{7}$ in one of their results. Unfortunately, their result came with some technical restrictions which meant that it could apparently not be combined with the GPY method to prove bounded gaps between primes\footnote{There is a technical discussion of exactly why not in the book \emph{Opera Cribro} by Friedlander and Iwaniec, pages 408--409.}. (As an aside, one of the main ingredients in these works of BFFI were certain complicated estimates for sums of Kloosterman sums due to Deshouillers and Iwaniec, and coming from the analytic theory of automorphic forms; some people apparently refer to this era as \emph{Kloostermania}.)

What Zhang succeeded in doing is modifying, in a quite nontrivial way, both the GPY method and the BFFI ideas so that they meet in the middle. 

As it turns out, an equivalent modification of the GPY method had already been published by Motohashi and Pintz. They observed that in the Elliott Halberstam conjecture $\EH(\theta)$ one does not need the primes to be nicely distributed modulo $d$ for most $d$ up to $X^{\theta}$. Rather one only needs\footnote{In fact one can get away with a still weaker property, in which one need only understand the number of primes congruent to $c$ mod $d$ for $c$ varying in a small set of residue classes modulo $d$ which varies in a multiplicative fashion with $d$.} this to be so for most \emph{smooth} values of $d$, that is to say values of $d$ with no prime factors bigger than $X^{\delta}$ for some very small $\delta$.  In particular, one does not need to know anything at all about the case when $d$ is prime or almost prime.  Thus Motohashi and Pintz, and independently Zhang, established a result of the following form.
\begin{theorem}[Motohashi-Pintz, Zhang]
Suppose that the primes are satisfactorily distributed modulo $d$ for most smooth values of $d$ up to $X^{\theta}$, for some value of $\theta > \frac{1}{2}$. Then we have bounded gaps between primes.
\end{theorem}
We could, if we wanted, call the fact that the primes are satisfactorily distributed modulo $d$ for most smooth values of $d$ up to $X^{\theta}$ the \emph{weak} Elliott-Halberstam conjecture, and denote it $\EH_{\weak}(\theta)$. To reiterate, then, Motohashi-Pintz-Zhang prove that if we have $\EH_{\weak}(\theta)$ for any $\theta > \frac{1}{2}$ then we still get bounded gaps between primes. 

Incidentally, now might be a good time to remark that the bound for the gap $H$ depends on how close $\theta$ is to $\frac{1}{2}$, the relation being very roughly of the form $H \sim (\theta - \frac{1}{2})^{-3/2}$ with a suitably optimised version of the argument, based on work of Conrey, Farkas, Pintz and Rev\'esz.

The heart of Zhang's advance is the following result.
\begin{theorem}[Zhang]\label{zhang-tech-1}
We have the weak Elliott-Halberstam conjecture $\EH_{\weak}(\theta)$ for $\theta = \frac{1}{2} + \frac{1}{1168}$. That is, the primes are satisfactorily distributed modulo $d$ for most smooth values of $d$ up to $X^{\theta}$. Hence, we have bounded gaps between primes.
\end{theorem}
We now turn to a few more details. First we say something about the GPY method (we shan't say anything here about its modification due to Motohashi-Pintz and Zhang). This uses ideas related to the \emph{Selberg sieve}. Then, we shall say a few words about the very technically demanding proof of Theorem \ref{zhang-tech-1}. The key words here are bilinear forms, Kloosterman sums and deep bounds of Bombieri and Birch coming from Deligne's proof of the Riemann hypothesis over finite fields. There is no use of \emph{automorphic form} bounds in Zhang's argument, and this is where he deviates somewhat from many of the papers of BFFI (though there are closely related ideas in a paper of Friedlander and Iwaniec on the divisor function).\vspace{11pt}

\textsc{The GPY method.} What has the distribution of primes modulo $d$ got to do with finding small gaps between primes? Exposing this hitherto unseen connection was the remarkable advance of Goldston, Pintz and Y{\i}ld{\i}r{\i}m. 

GPY in fact prove a result that is strictly stronger than bounded gaps between primes. We say that a $k$-tuple of integers $h_1,\dots,h_k$ is \emph{admissible} if there is no obvious ``congruence'' or ``local'' reason why $n + h_1, \dots, n+ h_k$ cannot all be prime for infinitely many $n$. For example, $\{h_1, h_2, h_3\} = \{0, 2, 4\}$ is not admissible, because at least one of these numbers is divisible by $3$, whereas $\{h_1, h_2, h_3\} = \{0,2,6\}$ is admissible (though no-one has the slightest idea how to \emph{prove} that $n, n+2, n+6$ are all prime for infinitely many $n$). A moment's thought convinces one that the natural criterion for admissibility is that, for each prime $p$, the set $\{h_1,\dots, h_k\}$ omits at least one residue class modulo $p$. If $n_*$ is this class then $n + h_1,\dots, n+ h_k$ could perhaps all be prime when $n = -n_* \md{p}$, because none of these numbers is divisible by $p$.

Here is the stronger statement that GPY proved.
\begin{theorem}\label{mpz-2}
Suppose we have the Elliott-Halberstam conjecture $\EH(\theta)$ for some $\theta > \frac{1}{2}$, that is to say the primes are nicely distributed modulo $d$ for most values of $d$ up to $X^{\theta}$. Suppose that $k \geq k_0(\theta)$ is sufficiently large. Then for any admissible $k$-tuple $h_1,\dots, h_k$, there are infinitely many $n$ for which at least two of $n+h_1,\dots, n+h_k$ are prime.
\end{theorem}
The modification of Motohashi-Pintz-Zhang is to show that this is still true if we instead assume just the weak Elliott-Halberstam conjecture $\EH_{\weak}(\theta)$, but we shall not be saying anyting further on the subject of this modification here.

To see why Theorem \ref{mpz-2} implies bounded gaps between primes, one need only note that there are admissible $k$-tuples for every $k$. Indeed the set of all prime numbers between $M$ and $2M$ will be admissible for all sufficiently large $M$, and the number of primes in this range grows without bound by the prime number theorem or in fact by weaker statements.

It should be pointed out that finding tightly packed admissible $k$ tuples, which is necessary to elucidate the relationship between $k$ and the prime gap $H$, brings one into contact with some thorny unsolved problems. To a large extent one must rely on computations to optimise this dependence for any particular value of $k$. 

We'll sketch the very broad outline of the proof of Theorem \ref{mpz-2}. It uses sieve theory, the branch of analytic number theory that has ultimately grown out of a serious study of the Sieve of Eratosthenes. Something that has been learned, rather painfully, over the last 100 years is that 
\[ \mbox{\emph{Almost primes are much easier to deal with than primes}}.\]
An $r$-almost prime is a product of at most $r$ primes. Fix an admissible tuple $\{h_1,\dots, h_k\}$. The rough idea of GPY is to choose an appropriate value of $r > k$ and try to compute the expected number of $n+h_1,\dots, n+ h_k$ that are prime, when $n$ is selected at random from those $n$ for which the product $(n+ h_1)\dots (n+ h_k)$ is an $r$-almost prime. If we can show that this is $> 1$ then, for some value of $n$, there must be two primes amongst the $n + h_i$.

To be a little more formal about this, write $\nu(n) = 1$ if $(n+h_1)\dots (n + h_k)$ is an almost prime and 0 otherwise. Let $X$ be an arbitrary large quantity. Then what we are interested in is the ratio
\begin{equation}\label{conditional} \frac{\sum_{X \leq n < 2X} (\Lambda(n + h_1) + \dots + \Lambda(n+ h_k)) \nu(n)}{\log X\sum_{X \leq n < 2X} \nu(n)}.\end{equation} (Recall that $\Lambda$ is basically the characteristic function of the primes weighted by $\log$.) If  we can show that this is $> 1$, we will then know that for some $n \in [X, 2X)$ at least two of $n+h_1,\dots, n+ h_k$ are prime, and this will conclude the proof of Theorem \ref{mpz-2}.

To elaborate this idea, one must be able to estimate the numerator and denominator of \eqref{conditional}. Now we come to a completely crucial idea, invented by Selberg. The idea is that there are weight functions $\nu(n)$ which behave morally rather like the characteristic function of the almost primes (or the set of $n$ for which $(n+h_1)\dots (n + h_k)$ is almost-prime), but which are much easier to compute with. 

Let $D$, $1 < D < X$, be a parameter and consider the function
\[ \nu(n) = \big(\sum_{\substack{d | n \\ d \leq D}} \lambda_d\big)^2,\] where at the moment $(\lambda_d)_{1 \leq d \leq D}$ is any set of real numbers with $\lambda_1 = 1$. The weight $\nu(n)$ is always nonnegative, and furthermore if $n$ is prime and between $D$ and $X$ then $\nu(n) = 1$. The reason is that in this case, the only divisors $d$ of $n$ are $1$ and $n$, and of these only $d = 1$ satisfies $d \leq D$. For this reason we'll call $\nu$ a \emph{majorant for the primes}.

Let's see how we can use a majorant like this to study a classical problem called the Brun-Titchmarsh problem: that of estimating from above the number of primes in a range $[X_0, X_0 + X)$. Now provided that\footnote{This should be read as ``a bit smaller than $X^{1/2}$''. In fact, one would require a condition like $D < X^{1/2 - \eps}$ for some $\eps > 0$.} $D = o(X^{1/2})$, we can compute an asymptotic for the average value of $\nu(n)$ over $X_0 \leq n < X_0 + X$ and hence get an upper bound for the number of primes in this range. To see why the condition $D = o(X^{1/2})$ is critical, we need to do an actual calculation (though a very short one):
\begin{align}\nonumber \sum_{X_0 \leq n <  X_0 + X} \nu(n) & = \sum_{X_0 < n \leq X_0 + X}\!\! \big(\sum_{\substack{d | n \\ d \leq D}} \lambda_d\big)^2 \\ & = \sum_{d,d' \leq D} \lambda_d \lambda_{d'}\!\!\!\! \sum_{\substack{X_0 \leq n < X_0 + X\\ d, d' | n}} 1.\label{actual-comp}\end{align}

Now the inner sum counts how many $n$ there are in the range $X_0 \leq n < X_0 + X$ for which both $d$ and $d'$ divide $n$, or equivalently for which the lowest common multiple $[d,d']$ divides $n$. Note that $[d,d'] \leq D^2$. Hence if $D = o(X^{1/2})$ then $[d,d'] = o(X)$ and the answer is essentially $X/[d,d']$. (Imagine you were asked how many multiples of 2014 there are in the interval $[10^{10}, 10^{10} + 10^5]$: since $\frac{10^5}{2014} = 49.65\dots$ it's either 49 or 50, but in either case $49.65$ is a good approximation.)

Therefore
\begin{equation}\label{eq4} \sum_{X_0 \leq n < X_0 + X} \nu(n) \approx \sum_{d,d' \leq D} \lambda_d \lambda_{d'} \frac{X}{[d,d']}.\end{equation}
If, however, $D > X^{1/2}$, then we could have (indeed we will often have) $[d,d'] > X$, so the answer might be 0, or it might be 1. (How many multiples of 2014 are there in the interval $[10^{10}, 10^{10} + 10^2]$?) We cannot say which without carefully inspecting $X_0$, and getting a usable expression is impossible without further ideas. Let us say that $X^{1/2}$ is the \emph{sieving limit} for this problem, and we call $D$ the \emph{sieving level}. Note that the larger we can take $D$, the more flexibility we have in choosing the weights $\lambda_d$. This ought to lead to the majorant $\nu$ being a better approximant to the characteristic function of the primes themselves. 

Miraculously, even though we are forced to take $D = o(X^{1/2})$ there are choices of the weights $\lambda_d$ for which $\nu$ is a reasonably good approximation to the characteristic function of the primes. We can find such $\lambda_d$ by minimising the quadratic form in \eqref{eq4} subject to $\lambda_1 = 1$, a routine exercise albeit one requiring some facility with M\"obius inversion. When this is done, we find that in fact, for this choice of weights, 
\[ \sum_{X_0 < n \leq X_0 + X} \nu(n) \approx \frac{2X}{\log X}.\]
When $X_0 = 0$ the majorant $\nu$ only overestimates the number of primes by a factor of $2$. If one looks very carefully at $\nu(n)$ then one sees that it resembles a sort of characteristic function of almost primes, though it is best not to pursue this line of thought too far, leaving it perhaps as motivation for calculating somewhat blindly with $\nu$. Note in particular that $\nu$ will not in general be $\{0,1\}$-valued. Note, by the way, that we have proven that the number of primes in $[X_0, X_0 + X]$ is at most about $2X/\log X$ for \emph{any} $X_0$, but that is another story.

It turns out that even if $D$ is a very small power of $X$ we still get a majorant $\nu$ that overestimates the primes by a constant factor, although this constant gets worse as $D$ becomes smaller. 

Returning to our main discussion, recall \eqref{conditional}. In the light of the crash course on the Selberg sieve we have just given, it is natural to consider defining
\begin{equation}\label{nu-def} \nu(n) = \big( \sum_{\substack{d | (n+h_1)\dots (n + h_k) \\ d \leq D}} \lambda_d \big)^2\end{equation} for appropriate weights $\lambda_d$, where $D$ is as big as possible.  To recap, one should think of $\nu$ as telling us the extent to which $(n+h_1) \dots (n+h_k)$ is almost prime, though to attach any more precise meaning to such a statement one must be more precise about the nature of the $\lambda_d$, about which I shall not say any more.

By a small modification of the above reasoning, one can compute the denominator in \eqref{conditional} provided that the sieving level $D$ is $o(X^{1/2})$. What, however, of the numerator? It may be split into terms of the form
\[ \sum_{n \leq X} \Lambda(n + h) \nu(n)\] for $h = h_1,\dots, h_k$.  Trying to repeat the computation in \eqref{actual-comp} above, we instead arrive at the expression
\begin{equation}\label{eq5}  \sum_{d_1,d_2 \leq D} \lambda_{d_1} \lambda_{d_2} \!\!\!\!\!\!\sum_{\substack{n \leq X \\ d_1,d_2 | (n+h_1)\dots (n + h_k)}} \Lambda(n).\end{equation}
To understand this sum, we need to know how the primes (weighted using the von Mangoldt function $\Lambda$) behave modulo $d = [d_1, d_2]$, and this quantity may be as large as $D^{2}$. If we know the Elliott-Halberstam conjecture $\EH(\theta)$ (primes are nicely distributed modulo $d$ for most $d$ up to $X^{\theta}$) then we will be fine so long as $D = o(X^{\theta/2})$. Unconditionally, that is to say using just the Bombieri-Vinogradov theorem, we may only take the sieving level $D$ to be about $X^{1/4}$, which means $\nu(n)$ gives a weaker notion of $(n+h_1) \dots (n + h_k)$ being almost prime. 

At this point one must dirty the hands by doing an actual computation of the numerator and denominator of \eqref{conditional} with a judicious choice of the weights $\lambda_d$. Making a sensible (by which we mean more-or-less optimal) choice, and taking $D$ to be almost $X^{\theta/2}$, one eventually computes that the ratio in \eqref{conditional} is 
\[ \frac{2\theta}{(1 + \frac{1}{\sqrt{k}})^2}.\]
Remember that $k$ is the number of elements in our admissible tuple $\{h_1,\dots, h_k\}$. I should say that I don't think I could motivate the result of this computation particularly well, if at all, even to an expert audience. I'm not such there even \emph{is} a conceptual explanation of it -- you just have to do it.

Recall that we wanted the ratio to be greater than 1: this would give us bounded gaps between primes. Even if $k$, the number of elements in our admissible tuple, is very large one does not achieve this if $\theta \leq \frac{1}{2}$. This is pretty unfortunate, since we can only proceed unconditionally when $\theta \leq \frac{1}{2}$. As soon as $\theta$ is even a tiny bit larger than $\frac{1}{2}$, however, the ratio will indeed be larger than $1$ provided that $k$ is big enough, and we will get bounded gaps between primes as discussed above.

The value $\theta = \frac{1}{2}$ is thus a crucial barrier for the GPY method: with $\theta < \frac{1}{2}$ one gets very little, whilst with $\theta > \frac{1}{2}$ one obtains bounded gaps between primes. \vspace{11pt}

This concludes our cursory discussion of the GPY method, which links bounded gaps between primes to the distribution of primes in progressions. We turn now to the other side of the story, in which the aim is to understand as much about the latter as possible.\vspace{11pt}

\textsc{Primes in arithmetic progression.} We turn now to a description of Zhang's major advance, the proof of Theorem \ref{zhang-tech-1}. Let us recall the statement.
\begin{theorem}[Zhang]
We have the weak Elliott-Halberstam conjecture $\EH_{\weak}(\theta)$ for $\theta = \frac{1}{2} + \frac{1}{1168}$. That is, the primes are satisfactorily distributed modulo $d$ for most smooth values of $d$ up to $X^{\theta}$. 
\end{theorem}
We did not say exactly what satisfactorily distributed means, but it basically means that we are interested in showing that\footnote{In fact, this only needs to be shown when $c$ belongs to a specific and fairly small set of residue classes varying in a multiplicative fashion with $d$, but we will not mention this point again.} if $(c,d) = 1$ then 
\[ \sum_{\substack{x \leq X \\ x \equiv c \mdsub{d}}} \Lambda(x) \approx \frac{X}{\phi(d)}\] for most smooth $d < X^{\theta}$. Remember that by smooth we meant that all prime factors of $d$ are at most $X^{\delta}$ for some very small $\delta$. The way in which this is used is that it allows us to suppose that $d$ is \emph{well-factorable}, which means that for any $Q, R$ with $QR \sim d$ we can find a factorisation $d = qr$ with $q \approx Q$ and $r \approx R$. Exactly \emph{how} this property is used is not very easy to explain properly, but bear with us. One can, however, immediately observe that if $q, r$ are coprime then the condition $x \equiv c \md{d}$ is equivalent to conditions on $x$ modulo $q$ and modulo $r$, by the Chinese remainder theorem, so there is a sense in which we are reducing the size of the moduli and thereby making the problem simpler.

At this point in the exposition it is slightly convenient to work with averages rather than sums, which we notate using the probabilistic $\E$ notation, though there is nothing random in our discussion. Our task is more-or-less equivalent to estimating an average over primes of the form
\[ \E_{x \leq X} \Lambda(x) \psi(x) = \frac{1}{X} \sum_{x  \leq X} \Lambda(x) \psi(x),\] where in this case 
\begin{equation}\label{psi-def} \psi(x) = 1_{x \equiv c \mdsub{d}} -  \frac{1}{\phi(d)}1_{(x, d) = 1} ,\end{equation}
 and the goal is to show that this average is appreciably smaller than the ``trivial'' bound of about $1$ which comes from the prime number theorem. 

To attack averages like this, we introduce a notion that is central to additive prime number theory: that of expanding in terms of \emph{bilinear forms}. Suppose that instead of the above average we were instead asked to estimate
\begin{equation}\label{to-convolve} \E_{x \leq X} (\alpha \ast \beta)(x) \psi(x),\end{equation} where $\alpha, \beta$ are arithmetic functions with $|\alpha(m)|, |\beta(n)| \leq 1$, and with $\alpha(m)$ supported on the range $m \sim M$ and $\beta(n)$ on the range $n \sim N$, where $MN = X$. Here, $\ast$ denotes Dirichlet convolution, that is to say
\[ (\alpha \ast \beta)(x) = \sum_{mn = x} \alpha(m) \beta(n).\]
The function $\psi$ is completely arbitrary for the purposes of this discussion, except we assume that $|\psi(x)| \leq 1$ for all $x$. We even let $\psi$ be complex-valued. 

The sum \eqref{to-convolve} can more-or-less be rewritten as 
\[ \E_{m \sim M}\E_{n \sim N} \alpha(m) \beta(n) \psi(mn).\]
Applying the Cauchy-Schwarz inequality, the square of this is bounded by
\[ \E_{m \sim M} |\E_{n \sim N} \beta(n) \psi(mn)|^2 =  \E_{n , n'\sim N} \beta(n) \beta(n') \E_{m \sim M} \psi(mn) \overline{\psi(mn')}.\]
Applying Cauchy-Schwarz a second time, the square of \emph{this} is bounded by 
\begin{equation}\label{gowers} \E_{m, m' \sim M} \E_{n, n' \sim N} \psi(mn)\overline{\psi(mn') \psi(m'n)} \psi(m'n').\end{equation}
The key thing to note here is that the unspecified functions $\alpha, \beta$ have completely disappeared, and we are left staring at an expression involving only $\psi$. Perhaps we might hope to estimate it, the aim being to improve substantially on the trivial bound of $1$. If we can do this, we have a kind of ``certificate'' which asserts that $\psi$ always gives nontrivial cancellation in averages such as \eqref{to-convolve}, no matter what $\alpha$ and $\beta$ are.

There are two rather obvious barriers to this observation being at all useful, and they are the following.
\begin{enumerate} \item We in fact wish to estimate the average $\E_{x \leq X} \Lambda(x) \psi(x)$, but we have said nothing about the extent to which $\Lambda$ can be expressed in terms of Dirichlet convolutions $\alpha \ast \beta$, nor even offered any motivation for why this should be expected.
\item We are interested in a specific function $\psi$, given by \eqref{psi-def}. Why should this $\psi$ allow us to provide a certificate by exhibiting nontrivial cancellation in \eqref{gowers}?
\end{enumerate}

With regard to (i), many readers will know that $\Lambda = \mu \ast \log$, where $\mu$ is the M\"obius function. However, this turns out not to help greatly in the above scheme. The reason is that in practice we will only be able to estimate expressions such as \eqref{gowers} for quite restricted ranges of $M$ and $N$, usually with $M$ and $N$ close in size, and with the decomposition $\mu \ast \log$ there is no opportunity to seriously restrict these ranges.

A successful technique depends on a remarkable identity of Linnik, or rather on a kind of truncated variant of it due to Heath-Brown, which we shall not state. The identity states that 
\[ \frac{\Lambda(n)}{\log n} = \sum_k \frac{(-1)^k}{k} \tau'_k(n),\] where $\tau'_k(n)$ is the number of ways to factor $n = n_1 \dots n_k$ with $n_i > 1$ for all $i$. If one knows the definition and very basic properties of the $\zeta$ function, the proof is just a couple of lines long: observe that
\[ \log (\zeta(s)) = \log (1 + (\zeta(s) - 1)) = \sum_{k = 1}^{\infty} (\zeta(s) - 1)^k,\]
and compare coefficients of $n^{-s}$ on both sides. The right hand side is pretty obviously the right-hand side of Linnik's identity. As for the left-hand side, its derivative is $\zeta'(s)/\zeta(s)$, which is well-known to be $-\sum_n \Lambda(n) n^{-s}$. Indeed, it is precisely this relation that links primes and the $\zeta$-function. Integrating with respect to $s$, we obtain Linnik's identity.

Note that $\tau'_k(n)$ is a Dirichlet convolution of $k$ copies of the function which equals $0$ at $1$ and is $1$ everywhere else. Chopping the domain of this function into various ranges, we can indeed write $\Lambda$ as a sum of a number (not too large) of convolutions $\alpha \ast \beta$, with $\alpha(m)$ supported where $m \sim M$ and $\beta(n)$ where $n \sim N$,  with  considerable flexibility in arranging the ranges $M$ and $N$ we need to worry about. The precise arrangement of these ranges is a rather technical matter, and suffice it to say that Zhang classifies them into three different types (plus a somewhat trivial type): these are called Type I, II and III. Actually, the Type III sums in fact involve certain 4-fold convolutions $\alpha_1 \ast \alpha_2 \ast \alpha_3 \ast \alpha_4$.

What about (ii), that is to say the issue of obtaining a ``certificate'' for $\psi$ which certifies that averages such as
\[ \E_{x \leq X} (\alpha \ast \beta)(x) \psi(x)\] 
exhibit cancellation? One may note that \eqref{gowers} is certainly not \emph{always} $o(1)$. Rather trivially, when $\psi$ is the constant function $1$ we get no cancellation. The same is true if $|\psi(x)| = 1$ and if $\psi$ is multiplicative in the sense that $\psi(mn) = \psi(m)\psi(n)$, as can be easily checked. (One consequence of this observation is that the whole scheme we have just outlined is somewhat unhelpful for showing that $\Lambda$ does not correlate with a single Dirichlet character $\chi$.) However, our function $\psi(x) = 1_{x \equiv c \mdsub{d}} -  \frac{1}{\phi(d)}1_{(x, d) = 1}$ does not obviously exhibit multiplicative behaviour and therefore one can hope to produce a certificate for this $\psi$, at least on average over $d$.

The reader will not be surprised to hear that the above discussion was an oversimplification, although it captures something of the key ideas. There are other types of ``certificate'' than \eqref{gowers}. The key tools that are brought to bear on estimating averages such as $\E_{x \leq X}(\alpha \ast \beta)(x) \psi(x)$ with our particular choice of $\psi$ are:
\begin{enumerate} 
\item Cauchy-Schwarz, similar to the above;
\item Fourier expansion, for example of $\psi$;
\item Shifting the range of summation, that is to say replacing the average $\E_{n \sim N} F(n)$ by $\E_{n \sim N} \E_{|k| \leq K} F(n+k)$, which should be roughly equal to it if $K \ll N$;
\item Certain changes of variable and substitutions;
\item Completion of sums, that is to say replacing an incomplete sum $\sum_{x \in I} f(x)$ by sums $\sum_{x \in \Z/d\Z} f(x) e_d(hx)$, where $I \subset \Z/d\Z$ is an interval, and $e_d(x) := e^{2\pi i x/d}$. 
\end{enumerate}

Considerable extra flexibility is available under the ``well-factorable'' assumption that $d = qr$; instead of averaging over $d \leq X^{\theta}$, one now averages over both $q$ and $r$, and this affords still more opportunity to vary the application of the above four ingredients.

The application of (i) to (v) above (in various combinations) throws up other sums that need to be estimated. The most interesting such case for Zhang occurs in his treatment of the so-called Type III sums, where expressions such as the sum
\begin{equation}\label{bombieri-birch} \sum_{n, n', \ell \mdsub{d}} e_d (\frac{c_1}{\ell n} + \frac{c_2}{(\ell + k)n'} + h_1n + h_2 n')\end{equation}  are relevant. Here, $h_1, h_2, k$ are parameters, and the sum over $n,n',\ell$ is restricted somewhat, in particular to $n, n', \ell, \ell + k \neq 0$ so that it makes sense.

It is very hard to explain how an expression like this comes up without going through the applications of each of the techniques (i) to (v) (which all occur here) in turn. Suffice it to say that the $k$ comes from shifting as in (iii), the $h_1, h_2$ come from (v) (followed by an application of Cauchy-Schwarz) and the instances of $c/n$ ultimately come from an initial Fourier expansion of $1_{x \equiv c(d)}$ where $x = mn$.

Now it turns out that \eqref{bombieri-birch} is in fact bounded by (essentially) $d^{3/2}$, apart from in certain degenerate cases. This is about as good an estimate as one should hope for, since it represents essentially square-root cancellation, as the number of things being summed over is $d^3$. This is by no means a trivial fact, and depends on cohomological ideas of Deligne, and in particular the Riemann hypothesis over finite fields. When $d = p$ is prime, the sum \eqref{bombieri-birch} may perhaps be written more suggestively, to an algebraic geometer, as
\[ \sum_{(x_1, x_2, x_3,x_4) \in V} e_p(x_1 + x_2 + x_3 + x_4),\] where $V \subset \F_p^4$ is a variety $\alpha x_1 x_2 + \beta x_3 x_4 = \gamma x_1 x_2 x_3 x_4$. 

In other cases (the analysis of the Type I and II sums) more standard sums such as Kloosterman sums
\[ \sum_{x_1 x_2 = c} e_p(x_1 + x_2)\] come up. This particular sum is bounded by $2\sqrt{p}$, a famous bound of Weil which does admit an elementary (in the technical sense) proof due to Stepanov. This is not the case with the Bombieri-Birch bound, which needs the whole of the algebro-geometric machinery. 

The great majority of these ideas can be found in the work of combinations of Bombieri, Fouvry, Friedlander and Iwaniec. The sad thing, however, is that this excellent and optimal bound of $d^{3/2}$ \emph{just} fails to cancel out some losses coming from other man{\oe}uvres, particularly the completion of sums manoeuvre (v) which is rather costly if the length of $I$ is much less than $d$. Instead of a valid certificate for $\psi$, one simply recovers, essentially, the trivial bound for \eqref{gowers} -- not, perhaps, the most spectacular use of the deep machinery.

The crucial new innovation of Zhang is to exploit the presence of the factorisation $d = qr$ (which, remember, can be selected quite flexibly). One might imagine that, by the Chinese remainder theorem, one obtains a product of sums similar to \eqref{bombieri-birch} modulo $q$ and modulo $r$, leading to a bound of $(qr)^{3/2}$ and no eventual gain. What Zhang miraculously finds, however, is that by making sure the shift in (iii) is by \emph{a multiple of $r$}, the sum modulo $r$ is not of Bombieri-Birch type \eqref{bombieri-birch}, but degenerates to something like 
\[ \sum_{\substack{s_1, s_2, n \mdsub{r} \\ s_1, s_2, n \neq 0}} e_r(\frac{s_1 - s_2}{s_1 s_2n}).\] This exhibits \emph{better than square root cancellation}, being of size essentially $r$, as one can easily check. (In fact, this is a slight simplification: Zhang actually obtains a Ramanujan sum in which the variable $n$ is constrained to be coprime to $r$, but the better-than-square-root cancellation still holds). 

Thus instead of $(qr)^{3/2}$ Zhang gets instead $q^{3/2} r$, and the factor of $r^{1/2}$ thus saved is crucial in making the argument work and establishing Zhang's remarkable theorem.\vspace{11pt}

\textsc{Polymath 8a.} Shortly after Zhang's paper came out, Terence Tao orchestrated a collaborative project to reduce the value of $H$ as far as possible from Zhang's value of 70000000, and more generally to understand all aspects of Zhang's work. One of the many great things about this project, in my view, was that without it there could well have been thousands of papers improving $H$ in various different ways. One big achievement of this project was to increase Zhang's $\frac{1}{1168}$ to $\frac{7}{300}$, that is to say to prove $\EH_{\weak}(\frac{1}{2} + \frac{3}{700})$, and to refine various other aspects of the argument, including the sieve theory and the computation determination of admissible tuples, so as to eventually reduce $H$ to $4680$.

Another significant achievement of the project was to remove the dependence on the deepest algebro-geometric results, although if this was one's concern then the value of $H$ could only be taken to be $14950$. At this point, the entire argument could reasonably be presented from first principles in an advanced graduate course. (Personally, I found it extraordinary that initially one could only get bounded gaps between primes using the full force of Deligne's machinery, and no finite bound without it.)\vspace{11pt}

\textsc{Maynard--Tao.} As \emph{Polymath 8a} was nearing completion, James Maynard and Terry Tao simultaneously made a dramatic advance which vastly simplifies the whole argument. Students wishing to study bound- ed gaps between primes may now read the 23-page paper of Maynard, freely available at \[ \mbox{ \texttt{http://arxiv.org/abs/1311.4600}.}\] In the Maynard--Tao argument, one still needs information about the distribution of primes in progressions, but things now work with any positive value of $\theta > 0$: the GPY barrier of $\theta = \frac{1}{2}$ turned out to be somewhat illusory. The crucial new idea of Maynard and Tao is to consider a different weight function $\nu$. Instead of 
\[\nu(n) = \big( \sum_{\substack{d | (n+h_1)\dots (n + h_k) \\ d \leq D}} \lambda_d \big)^2,\] they consider the more general function
\[ \nu(n) = \big( \sum_{\substack{d_1 | (n+h_1), d_2 |(n + h_2), \dots d_k| (n + h_k) \\ d_1\cdots d_k \leq D}} \lambda_{d_1,\dots, d_k} \big)^2.\]
The task is again to choose the weights $\lambda_{d_1,\dots, d_k}$ so as to try and make the ratio \eqref{conditional} bigger than $1$. It turns out that with an appropriate choice the ratio can be made to behave roughly like $\frac{1}{2}\theta \log k$. For any fixed value of $\theta$ (and not just for $\theta > \frac{1}{2}$!) this can be made greater than $1$ by choosing $k$ big enough. In fact, it can be made larger than any preassigned quantity, and so the Maynard--Tao argument shows that there are infinitely many triples, quadruples, quintuples, ... of primes in bounded gaps. No such result was previously known. Maynard's paper also reduces the gap $H$ to 600. \vspace{11pt}

\textsc{Polymath 8b.} A second Polymath project, again led by Terence Tao, has been working on the problem in the light of the Maynard--Tao development. When the first version of this paper was submitted on 20/2/14 the record value of $H$ was 264, but when I came to correct some typos a few days later this value had already been reduced to $H = 252$. A significant achievement of this second project has been to show that if one assumes the \emph{full} Elliott-Halberstam conjecture $\EH(\theta)$ for every $\theta < 1$ then one obtains $H = 6$, which is known to be optimal using anything like these methods (Maynard had earlier obtained $H = 12$ on the same assumption, and GPY the value $H = 16$). Whether the Elliott-Halberstam conjecture is an easier target than the twin prime conjecture is a matter for debate, but at least one may say that the dimension of the space of conjectures has been somewhat reduced.

\end{document}